\documentclass[12pt]{article}%
\usepackage[intlimits]{amsmath}
\usepackage{amssymb}
\usepackage{latexsym}
\usepackage[T1]{fontenc}
\usepackage[sc]{mathpazo}
\usepackage{color}
\usepackage[colorlinks]{hyperref}
\usepackage{amsfonts}
\usepackage{microtype}
\usepackage{graphicx}%
\definecolor {refcol}{RGB}{40,0,255}
\hypersetup{colorlinks=true,allcolors=refcol}
\makeatletter
\newfont{\footsc}{cmcsc10 at 8truept}
\newfont{\footbf}{cmbx10 at 8truept}
\newfont{\footrm}{cmr10 at 10truept}
\newtheorem{theorem}{Theorem}

\newtheorem{corollary}[theorem]{Corollary}

\newtheorem{problem}[theorem]{Problem}
\newtheorem{proposition}[theorem]{Proposition}

\newenvironment{proof}[1][Proof]{\noindent{\textbf {#1}  }}  {\hfill$\Box$\bigskip}

\begin{document}

\title{\textbf{Spectral radius and Hamiltonicity of graphs with large minimum degree}}
\author{V. Nikiforov\thanks{Department of Mathematical Sciences, University of
Memphis, Memphis TN 38152, USA}\medskip\\{\small Dedicated to the memory of Miroslav Fiedler}}
\date{}
\maketitle

\begin{abstract}
\ This paper presents sufficient conditions for Hamiltonian paths and cycles
in graphs. Letting $\lambda\left(  G\right)  $ denote the spectral radius of
the adjacency matrix of a graph  $G,$ the main results of the paper
are:\medskip

(1) Let $k\geq1,$ $n\geq k^{3}/2+k+4,$ and let $G$ be a graph of order $n$,
with minimum degree $\delta\left(  G\right)  \geq k.$ If
\[
\lambda\left(  G\right)  \geq n-k-1,
\]
then $G$ has a Hamiltonian cycle, unless $G=K_{1}\vee(K_{n-k-1}+K_{k})$ or
$G=K_{k}\vee(K_{n-2k}+\overline{K}_{k}).$\medskip

(2) Let $k\geq1,$ $n\geq k^{3}/2+k^{2}/2+k+5,$ and let $G$ be a graph of order
$n$, with minimum degree $\delta\left(  G\right)  \geq k.$ If
\[
\lambda\left(  G\right)  \geq n-k-2,
\]
then $G$ has a Hamiltonian path, unless $G=K_{k}\vee(K_{n-2k-1}+\overline
{K}_{k+1})$ or $G=K_{n-k-1}+K_{k+1}$\medskip

In addition, it is shown that in the above statements, the bounds on $n$ are
tight within an additive term not exceeding $2$.\medskip

\medskip


\end{abstract}

\section{Introduction}

In this paper we present sufficient spectral conditions for Hamiltonian paths
and cycles in graphs with large minimum degree.

Write $\lambda\left(  G\right)  $ for the spectral radius of the adjacency
matrix of a graph $G$. In 2010, Fiedler and Nikiforov \cite{FiNi10} gave some
bounds on $\lambda\left(  G\right)  $ implying the existence of Hamiltonian
paths and cycles in $G$. This work motivated further research, as could be
seen, e.g., in \cite{Ben15, GeNi16,Li15, LiNi15, LLT12, LSX15, NiGe15, Zho10}.

The aim of this paper is to extend some recent results by Benediktovich
\cite{Ben15}, Li and Ning \cite{LiNi15}, and Ning and Ge \cite{NiGe15}. To
state those results we need to introduce three families of extremal graphs,
which are denoted by $L_{k}\left(  n\right)  ,$ $M_{k}\left(  n\right)  $ and
$N_{k}\left(  n\right)  $.

Write $K_{s}$ and $\overline{K}_{s}$ for the complete and the edgeless graphs
of order $s.$ Given graphs $G$ and $H,$ write $G\vee H$ for their join and
$G+H$ for their disjoint union.\medskip\

\textbf{The graphs }$L_{k}\left(  n\right)  .$ For any $k\geq1$ and
$n\geq2k+1,$ let%
\[
L_{k}\left(  n\right)  :=K_{1}\vee\left(  K_{n-k-1}+K_{k}\right)  .
\]
Thus, the graph $L_{k}\left(  n\right)  $ consists of a $K_{n-k}$\ and a
$K_{k+1}$ sharing a single vertex.\medskip

\textbf{The graphs}\emph{ }$M_{k}\left(  n\right)  .$ For any $k\geq1$ and
$n\geq2k+1,$ let
\[
M_{k}\left(  n\right)  :=K_{k}\vee\left(  K_{n-2k}+\overline{K}_{k}\right)  .
\]
Thus, the graph $M_{k}\left(  n\right)  $ consists of a $K_{n-k}$\ and a set
of $k$ independent vertices all joined to some $k$ vertices of the $K_{n-k}%
$.\medskip

\textbf{The graphs}\emph{ }$N_{k}\left(  n\right)  .$ For any $k\geq1$ and
$n\geq2k+1,$ let
\[
N_{k}\left(  n\right)  :=K_{k}\vee\left(  K_{n-2k-1}+\overline{K}%
_{k+1}\right)  .
\]
Thus, the graph $N_{k}\left(  n\right)  $ consists of a $K_{n-k-1}$\ and a set
of $k+1$ independent vertices all joined to some $k$ vertices of the
$K_{n-k-1}.\medskip$

Note that for any admissible $k$ and $n,$ the graphs $L_{k}\left(  n\right)  $
and $M_{k}\left(  n\right)  $ contain no Hamiltonian cycle and $N_{k}\left(
n\right)  $ contains no Hamiltonian path, whereas the minimum degree of each
of them is exactly $k.\medskip$

The graphs $M_{k}\left(  n\right)  $ and $N_{k}\left(  n\right)  $ were used
by Erd\H{o}s \cite{Erd62} as extremal graphs in his results on Hamiltonicity
of graphs with large minimum degree. Moreover, recently Li and Ning
\cite{LiNi15} showed that $M_{k}\left(  n\right)  $ and $N_{k}\left(
n\right)  $ are relevant also for some spectral analogs of Erd\H{o}s's results:

\begin{theorem}
[Li, Ning \cite{LiNi15}]\label{tln}Let $k\geq0$ and $G$ be a graph of order
$n,$ with minimum degree $\delta\left(  G\right)  \geq k.$

(1) If $n\geq\max\left\{  6k+10,(k^{2}+7k+8)/2\right\}  $ and
\[
\lambda\left(  G\right)  \geq\lambda\left(  N_{k}\left(  n\right)  \right)  ,
\]
then $G$ has a Hamiltonian path, unless $G=N_{k}\left(  n\right)  $;

(2) If $k\geq1,$ $n\geq\max\left\{  6k+5,(k^{2}+6k+4)/2\right\}  ,$ and
\[
\lambda\left(  G\right)  \geq\lambda\left(  M_{k}\left(  n\right)  \right)  ,
\]
then $G$ has a Hamiltonian cycle, unless $G=M_{k}\left(  n\right)  .$
\end{theorem}

Theorem \ref{tln} seems as good as one can get, yet somewhat subtler and
stronger statements have been proved for $k=1,2.$

\begin{theorem}
[Ning, Ge \cite{NiGe15}]\label{tng}Let $n\geq7$ and $G$ be a graph of order
$n$, with minimum degree $\delta\left(  G\right)  \geq1.$ If
\[
\lambda\left(  G\right)  >n-3,
\]
then $G$ has a Hamiltonian path, unless $G=N_{1}\left(  n\right)  .$
\end{theorem}

\begin{theorem}
[Benediktovich \cite{Ben15}]\label{tb}Let $n\geq10$ and let $G$ be a graph of
order $n$, with minimum degree $\delta\left(  G\right)  \geq2.$ If
\[
\lambda\left(  G\right)  \geq n-3,
\]
then $G$ has a Hamiltonian cycle, unless $G=L_{2}\left(  n\right)  $ or
$G=M_{2}\left(  n\right)  .$
\end{theorem}

Note that in our renditions of Theorems \ref{tng} and \ref{tb}, a few details
have been suppressed from the original theorems in order to highlight them as
instances of more general statements, in which $\delta\left(  G\right)  $ is
bounded by a parameter.\medskip

Thus, here we propose the following two theorems, which generalize Theorems
\ref{tng} and \ref{tb}, and strengthen Theorem \ref{tln} for $n$ sufficiently large:

\begin{theorem}
\label{mtc}Let $k\geq1,$ $n\geq k^{3}/2+k+4,$ and let $G$ be a graph of order
$n$, with minimum degree $\delta\left(  G\right)  \geq k.$ If
\[
\lambda\left(  G\right)  \geq n-k-1,
\]
then $G$ has a Hamiltonian cycle, unless $G=L_{k}\left(  n\right)  $ or
$G=M_{k}\left(  n\right)  .$
\end{theorem}

\begin{theorem}
\label{mtp}Let $k\geq1,$ $n\geq k^{3}/2+k^{2}/2+k+5,$ and let $G$ be a graph
of order $n$, with minimum degree $\delta\left(  G\right)  \geq k.$ If
\[
\lambda\left(  G\right)  \geq n-k-2,
\]
then $G$ has a Hamiltonian path, unless $G=N_{k}\left(  n\right)  $ or
$G=K_{n-k-1}+K_{k+1}.$
\end{theorem}

We shall give independent, self-contained proofs of Theorems \ref{mtc} and
\ref{mtp}, although\ many smart ideas could be readily borrowed from each of
the papers \cite{Ben15}, \cite{LiNi15}, and \cite{NiGe15}. Crucial points of
our arguments are based on the following straightforward theorems, whose
proofs are nevertheless long and technical:

\begin{theorem}
\label{t2}Let $k\geq1,$ $n\geq k^{3}/2+k+4,$ and let $G$ be a graph of order
$n$, with minimum degree $\delta\left(  G\right)  \geq k.$

(i) If $G$ is a subgraph of $M_{k}\left(  n\right)  $, then $\lambda\left(
G\right)  <n-k-1$, unless $G=M_{k}\left(  n\right)  .$

(ii) If $G$ is a subgraph of $L_{k}\left(  n\right)  $, then $\lambda\left(
G\right)  <n-k-1$, unless $G=L_{k}\left(  n\right)  .$
\end{theorem}

\begin{theorem}
\label{t3}Let $k\geq1,$ $n\geq k^{3}/2+k^{2}/2+k+5,$ and let $G$ be a graph of
order $n$, with minimum degree $\delta\left(  G\right)  \geq k.$

(i) If $G$ is a subgraph of $N_{k}\left(  n\right)  $, then $\lambda\left(
G\right)  <n-k-2$, unless $G=N_{k}\left(  n\right)  .$

(ii) If $G$ is a subgraph of $K_{n-k-1}+K_{k+1}$, then $\lambda\left(
G\right)  <n-k-2$, unless $G=K_{n-k-1}+K_{k+1}.$
\end{theorem}

One can ask how tight are the lower bounds on $n$ in the premises of Theorems
\ref{mtc}--\ref{t3}. This question has been brought up by Ge and Ning in
\cite{GeNi16}, after they reduced by half the bound on $n$ in a previous
version of Theorem \ref{t2}. In fact, the bounds turn out to be tight within
an additive term not exceeding $2$:

\begin{proposition}
\label{pro1} If $k\geq2$ and $2k+1\leq n\leq k^{3}/2+k+1$, there exists a
graph $G$ of order $n$, with minimum degree $\delta\left(  G\right)  \geq k$
such that $G$ is a proper subgraph of $M_{k}\left(  n\right)  $, and
$\lambda\left(  G\right)  >n-k-1.$
\end{proposition}

\begin{proposition}
\label{pro2} If $k\geq2$ and $2k+1\leq n\leq k^{3}/2+k^{2}/2+k+2$, there
exists a graph $G$ of order $n$, with minimum degree $\delta\left(  G\right)
\geq k$ such that $G$ is a proper subgraph of $N_{k}\left(  n\right)  $, and
$\lambda\left(  G\right)  >n-k-1.$
\end{proposition}

Clearly, Propositions \ref{pro1} and \ref{pro2} leave room for some ultimate
nitpicking, which we skip.$\medskip$

The rest of the paper is structured as follows: In Section \ref{s2}, we
introduce some notation, recall some details about graph closure, and state a
few results that will be used in the proofs of\ Theorems \ref{mtc}-\ref{t3}
and Propositions \ref{pro1} and \ref{pro2}. The proofs themselves are given in
Section \ref{s3}. The last Section \ref{s4} is dedicated to a brief discussion
and some open problems.

\section{\label{s2}Notation and preliminaries}

For graph notation and terminology undefined here we refer the reader to
\cite{Bol98}. We write $A\left(  G\right)  $ for the adjacency matrix of a
graph $G$, and denote the quadratic form of $A\left(  G\right)  $ by
$\left\langle A\left(  G\right)  \mathbf{x},\mathbf{x}\right\rangle ,$ where
$\mathbf{x}$ is a vector of size equal to the order of $G.$ Note that if $G$
is of order $n$ and $\mathbf{x}:=\left(  x_{1},\ldots,x_{n}\right)  ,$ then
\[
\left\langle A\left(  G\right)  \mathbf{x},\mathbf{x}\right\rangle
=2\sum_{\left\{  i,j\right\}  \in E\left(  G\right)  }x_{i}x_{j}.
\]

If $G$ is a graph of order $n,$ we write $d_{1},\ldots,d_{n}$ for the degrees
of $G$ indexed in ascending order$.$

A graph $G$ is called \emph{Hamiltonian-connected} if for any two vertices $u$
and $v$ of $G,$ there is a Hamiltonian path in $G$ whose ends are $u$ and
$v.\medskip$

We shall need the concept of graph closure, used implicitly by Ore in
\cite{Ore61, Ore63}, and developed further by Bondy and Chv\'{a}tal in
\cite{BoCh76}: Fix an integer $k\geq0.$ Given a graph $G,$ perform the
following operation: if there are two nonadjacent vertices $u$ and $v$ with
$d_{G}\left(  u\right)  +d_{G}\left(  v\right)  \geq k,$ add the edge $uv$ to
$E\left(  G\right)  .$ A $k$\emph{-closure} of $G$ is a graph obtained from
$G$ by successively applying this operation as long as possible. As it turns
out, the $k$-closure of $G$ is unique, that is to say, it does not depend on
the order in which edges are added; see \cite{BoCh76} for details.\medskip

Write $\mathrm{cl}_{k}\left(  G\right)  $ for the $k$-closure of $G$ and note
its main property:\medskip

\emph{If }$u$\emph{ and }$v$\emph{ are nonadjacent vertices of }%
$\mathrm{cl}_{k}\left(  G\right)  ,$\emph{ then }$d_{\mathrm{cl}_{k}\left(
G\right)  }\left(  u\right)  +d_{\mathrm{cl}_{k}\left(  G\right)  }\left(
v\right)  \leq k-1.$\medskip

The usefulness of graph closure is demonstrated by the following facts, due
essentially to Ore \cite{Ore61, Ore63}:

- A graph $G$ of order $n$ has a Hamiltonian path if and only if
$\mathrm{cl}_{n-1}\left(  G\right)  $ has one.

- A graph $G$ of order $n$ has a Hamiltonian cycle if and only if
$\mathrm{cl}_{n}\left(  G\right)  $ has one.

- A $2$-connected graph $G$ of order $n$ is Hamiltonian-connected if and only
if $\mathrm{cl}_{n+1}\left(  G\right)  $ is Hamiltonian-connected.\medskip

For convenience we restate the last two statements in more usable form.

\begin{theorem}
[Ore \cite{Ore63}]\label{ot1}If $G$ is a graph of order $n$ and $d_{u}%
+d_{v}\geq n$ for any two distinct nonadjacent vertices $u$ and $v$, then $G$
has a Hamiltonian cycle.
\end{theorem}

\begin{theorem}
[Ore \cite{Ore63}]\label{ot}If $G$ is a $2$-connected graph of order $n$ and
$d_{u}+d_{v}\geq n+1$ for any two distinct nonadjacent vertices $u$ and $v$,
then $G$ is Hamiltonian-connected.
\end{theorem}

We shall also need two classical results of Chv\'{a}tal \cite{Chv72} on
Hamiltonicity of graphs.\

\begin{theorem}
[Chv\'{a}tal \cite{Chv72}]\label{cht}Let $G$ be a graph of order $n,$ with
degrees $d_{1},\ldots,d_{n}$. If $G$ has no Hamiltonian cycle, then there is
an integer $s<n/2$ such that $d_{s}\leq s$ and $d_{n-s}\leq n-s-1.$
\end{theorem}

\begin{corollary}
[Chv\'{a}tal \cite{Chv72}]\label{chc}Let $G$ be a graph of order $n,$ with
degrees $d_{1},\ldots,d_{n}$. If $G$ has no Hamiltonian path, then there is an
integer $s<\left(  n+1\right)  /2$ such that $d_{s}\leq s-1$ and
$d_{n-s+1}\leq n-s-1.$
\end{corollary}

Finally, we shall need the following inequality, proved in \cite{Nik02}:

\begin{theorem}
[\cite{HSF01}, \cite{Nik02}]\label{tin}If $G$ is a graph of order $n$, with
$m$ edges, and minimum degree $\delta,$ then
\begin{equation}
\lambda\left(  G\right)  \leq\frac{\delta-1}{2}+\sqrt{2m-n\delta+\frac{\left(
\delta+1\right)  ^{2}}{4}}. \label{in}%
\end{equation}

\end{theorem}

For connected graphs inequality (\ref{in}) has been proved independently by
Hong, Shu and Fang\ in \cite{HSF01}.

The following observation is often useful in applications of inequality
(\ref{in}).

\begin{proposition}
[\cite{HSF01}, \cite{Nik02}]\label{pro}If $2m\leq n\left(  n-1\right)  $, the
function
\[
f\left(  x\right)  =\frac{x-1}{2}+\sqrt{2m-nx+\frac{\left(  x+1\right)  ^{2}%
}{4}}%
\]
is decreasing in $x$ for $x\leq n-1.$
\end{proposition}

\section{\label{s3}Proofs}

\subsection{Proof of Theorem \ref{t2}}

\begin{proof}
Set for short $\lambda:=\lambda\left(  G\right)  ,$ and let $\mathbf{x}%
=\left(  x_{1},\ldots,x_{n}\right)  $ be a positive unit eigenvector vector to
$\lambda.$ Recall that Rayleigh's principle implies that
\[
\lambda=\left\langle A\left(  G\right)  \mathbf{x},\mathbf{x}\right\rangle .
\]
The idea of the proofs below exploits the fact that both $M_{k}\left(
n\right)  $ and $L_{k}\left(  n\right)  $ consist of a $K_{n-k},$ together
with an \textquotedblleft outgrowth\textquotedblright\ of bounded order. It
turns out that if $n$ is large, the total contribution to $\left\langle
A\left(  G\right)  \mathbf{x},\mathbf{x}\right\rangle $ of all edges incident
to the \textquotedblleft outgrowth\textquotedblright\ is less than the
contribution of a single edge of the $K_{n-k}.$

Now we give the details.$\medskip$

\emph{Proof of (i)}

Assume that $G$ is a proper subgraph of $M_{k}\left(  n\right)  .$ Clearly, we
may assume that $G$ is obtained by omitting just one edge $\left\{
u,v\right\}  $ of $M_{k}\left(  n\right)  .$

Write $X$ for the set of vertices of $M_{k}\left(  n\right)  $ of degree $k,$
let $Y$ be the set of their neighbors, and let $Z$ be the set of the remaining
$n-2k$ vertices of $M_{k}\left(  n\right)  .$

Since $\delta\left(  G\right)  \geq k,$ we see that $G$ must contain all the
edges between $X$ and $Y.$ Therefore, $\left\{  u,v\right\}  \subset Y\cup Z,$
with three possible cases: (a) $\left\{  u,v\right\}  \subset Y$; (b) $u\in
Y,v\in Z$; (c) $\left\{  u,v\right\}  \subset Z.$ We shall show that case (c)
yields a graph of no smaller spectral radius than case (b), and that case (b)
yields a graph of no smaller spectral radius than case (a).

Indeed, by symmetry, we have $x_{i}=x_{j}$ for any $i,j\in X$; likewise,
$x_{i}=x_{j}$ for any $i,j\in Y\backslash\left\{  u,v\right\}  $ and for any
$i,j\in Z\backslash\left\{  u,v\right\}  .$ Thus, let
\begin{align*}
x &  :=x_{i},\text{ }i\in X,\\
y &  :=x_{i},\text{ }i\in Y\backslash\left\{  u,v\right\}  ,\\
z &  :=x_{i},\text{ }i\in Z\backslash\left\{  u,v\right\}  .
\end{align*}

Suppose that case (a) holds, that is, $\left\{  u,v\right\}  \subset Y.$
Choose a vertex $w\in Z,$ remove the edge $\left\{  v,w\right\}  $ and add the
edge $\left\{  u,v\right\}  .$ If $x_{w}>x_{v},$ swap the entries $x_{v}$ and
$x_{w}$; write $\mathbf{x}^{\prime}$ for the resulting vector.

First, note that $\mathbf{x}^{\prime}$ is a unit vector and that the obtained
graph $G^{\prime}$ is covered by case (b). \ We see that
\[
\left\langle A\left(  G^{\prime}\right)  \mathbf{x}^{\prime},\mathbf{x}%
^{\prime}\right\rangle -\left\langle A\left(  G\right)  \mathbf{x}%
,\mathbf{x}\right\rangle =\left(  x_{v}^{\prime}-x_{v}\right)  \sum_{i\in
X}x_{i}\geq0,
\]
and, by the Rayleigh principle, $\lambda\left(  G^{\prime}\right)  \geq
\lambda\left(  G\right)  ,$ as claimed$.$

Essentially the same argument proves that case (c) yields a graph of no
smaller spectral radius than case (b). Therefore, we may assume that $\left\{
u,v\right\}  \subset Z.$ Hence, the vertices $u$ and $v$ are symmetric, and so
$x_{u}=x_{v}.$ Set $t:=x_{u}$ and note that the $n$ eigenequations of $G$ are
reduced to four equations involving just the unknowns $x,$ $y,$ $z,$ and $t:$%
\begin{align}
\lambda x  &  =ky,\label{e1}\\
\lambda y  &  =kx+\left(  k-1\right)  y+\left(  n-2k-2\right)  z+2t,
\label{e2}\\
\lambda z  &  =ky+\left(  n-2k-3\right)  z+2t,\label{e3}\\
\lambda t  &  =ky+\left(  n-2k-2\right)  z. \label{e4}%
\end{align}
We find that
\begin{align}
x  &  =\frac{k}{\lambda}y,\label{xs}\\
z  &  =\left(  1-\frac{k^{2}}{\lambda\left(  \lambda+1\right)  }\right)
y,\nonumber\\
t  &  =\frac{\lambda+1}{\lambda+2}\left(  1-\frac{k^{2}}{\lambda\left(
\lambda+1\right)  }\right)  y. \label{etz}%
\end{align}

Further, note that if we remove all edges between $X$ and $Y$ and add the edge
$\left\{  u,v\right\}  $ to $G,$ we obtain the graph $K_{n-k}+\overline{K}%
_{k}.$ Letting $\mathbf{x}^{\prime\prime}$ be the restriction of $\mathbf{x}$
to $K_{n-k},$ we find that
\[
\left\langle A\left(  K_{n-k}\right)  \mathbf{x}^{\prime\prime},\mathbf{x}%
^{\prime\prime}\right\rangle =\left\langle A\left(  G\right)  \mathbf{x}%
,\mathbf{x}\right\rangle +2t^{2}-2k^{2}xy=\lambda+2t^{2}-2k^{2}xy.
\]
But since $\left\Vert \mathbf{x}^{\prime\prime}\right\Vert ^{2}=1-kx^{2}$, we
see that\footnote{This equation was suggested by Ge and Ning in \cite{GeNi16}
instead of the inequality $\left\Vert \mathbf{x}^{\prime\prime}\right\Vert <1$
used in a previous version of the present proof. }
\[
\lambda+2t^{2}-2k^{2}xy=\left\langle A\left(  K_{n-k}\right)  \mathbf{x}%
^{\prime\prime},\mathbf{x}^{\prime\prime}\right\rangle \leq\lambda\left(
K_{n-k}\right)  \left\Vert \mathbf{x}^{\prime\prime}\right\Vert ^{2}=\left(
n-k-1\right)  (1-kx^{2}).
\]
Assume for a contradiction that $\lambda\geq n-k-1.$ This assumption, together
with above inequality, yields
\[
\lambda+2t^{2}-2k^{2}xy\leq(1-kx^{2})\lambda,
\]
and therefore
\[
2t^{2}-2k^{2}xy\leq-kx^{2}\lambda.
\]
Now, (\ref{xs}) and (\ref{etz}) imply that
\[
\frac{k^{3}}{2\lambda}y^{2}\geq2t^{2}=\left(  \frac{\lambda+1}{\lambda
+2}\right)  ^{2}\left(  1-\frac{k^{2}}{\lambda\left(  \lambda+1\right)
}\right)  ^{2}y^{2}.
\]
Cancelling $y^{2}$ and applying Bernoulli's inequality to the right side, we
get
\begin{align*}
k^{3} &  \geq2\lambda\left(  1-\frac{1}{\lambda+2}\right)  ^{2}\left(
1-\frac{k^{2}}{\lambda\left(  \lambda+1\right)  }\right)  ^{2}>2\lambda
-\frac{4\lambda}{\lambda+2}-\frac{4k^{2}}{\lambda+1}\\
&  >2\lambda-\frac{4\lambda+4k^{2}}{\lambda+1}.
\end{align*}
Using the inequalities $\lambda\geq n-k-1\geq k^{3}/2+3,$ we easily find that
\[
\frac{4\lambda+4k^{2}}{\lambda+1}\leq6,
\]
and so,
\[
k^{3}>2\lambda-\frac{4\lambda+4k^{2}}{\lambda+1}>k^{3}+6-6=k^{3},
\]
a contradiction, completing the proof of (i).\medskip

\emph{Proof of (ii)}

As in (i), assume that $G$ is a subgraph of $L_{k}\left(  n\right)  $ obtained
by omitting just one edge $\left\{  u,v\right\}  $ of $L_{k}\left(  n\right)
.$ Recall that $L_{k}\left(  n\right)  $ consists of a $K_{n-k}$ and a
$K_{k+1}$ sharing a single vertex, say $w$. Let $Y$ be the set $\left\{
w\right\}  ,$ write $X$ for the set of vertices of $K_{k+1}$ that are distinct
from $w,$ and write $Z$ for the set of vertices of $K_{n-k}$ that are distinct
from $w.$

Clearly, the condition $\delta\left(  G\right)  \geq k$ implies that $\left\{
u,v\right\}  \subset Y\cup Z$; among the three possible placements of
$\left\{  u,v\right\}  ,$ the case $\left\{  u,v\right\}  \subset Z$ yields a
graph with maximum spectral radius$,$ so we assume that $\left\{  u,v\right\}
\subset Z.$ Now, by symmetry, we see that $x_{i}=x_{j}$ for any $i,j\in X$;
likewise, $x_{u}=x_{v}$ and $x_{i}=x_{j}$ for any $i,j\in Z\backslash\left\{
u,v\right\}  .$ Thus, let
\begin{align*}
x &  :=x_{i},\text{ }i\in X,\\
y &  :=x_{w},\text{ }i\in Y,\\
z &  :=x_{i},\text{ }i\in Z\backslash\left\{  u,v\right\}  ,\\
t &  :=x_{u}=x_{v}.
\end{align*}
The $n$ eigenequations of $G$ now reduce to the four equations%
\begin{align*}
\lambda x &  =y+\left(  k-1\right)  x,\\
\lambda y &  =kx+\left(  n-2k-3\right)  z+2t,\\
\lambda z &  =y+\left(  n-2k-4\right)  z+2t,\\
\lambda t &  =y+\left(  n-2k-3\right)  z.
\end{align*}
Hence, we find that
\begin{align*}
x &  =\frac{1}{\lambda-k+1}y,\\
z &  =\left(  1-\frac{k}{\left(  \lambda-k+1\right)  \left(  \lambda+1\right)
}\right)  y,\\
t &  =\frac{\lambda+1}{\lambda+2}\left(  1-\frac{k}{\left(  \lambda
-k+1\right)  \left(  \lambda+1\right)  }\right)  y.
\end{align*}

Further, if we delete all edges incident to vertices in $X$ and add the edge
$\left\{  u,v\right\}  ,$ we obtain the graph $K_{n-k}+\overline{K}_{k}.$

Assume for a contradiction that $\lambda\geq n-k-1$. Reasoning as in the proof
of (i), we get the inequality
\[
\frac{k\left(  k-1\right)  }{2\left(  \lambda-k+1\right)  ^{2}}y^{2}+\frac
{k}{\lambda-k+1}y^{2}>t^{2}=\left(  \frac{\lambda+1}{\lambda+2}\right)
^{2}\left(  1-\frac{k}{\left(  \lambda-k+1\right)  \left(  \lambda+1\right)
}\right)  ^{2}y^{2},
\]
which in turn yields%
\[
\frac{k\left(  k-1\right)  }{2\left(  \lambda-k+1\right)  }+k>\lambda
-k+1-\frac{2\left(  \lambda-k+1\right)  }{\lambda+2}-\frac{2k}{\lambda
+1}>\lambda-k-1.(1-kx^{2})
\]
Hence, using the inequality $\lambda\geq n-k-1\geq k^{3}/2+3,$ we find that
\begin{align*}
k\left(  k-1\right)   &  >2\left(  \lambda-k+1\right)  \left(  \lambda
-2k-1\right)  \\
&  >\left(  k^{3}-2k+8\right)  \left(  \frac{k^{3}}{2}-2k+2\right)  \\
&  >2k^{3}-4k+16.
\end{align*}
It is not hard to see that this inequality is a contradiction for $k\geq2,$
completing the proof of Theorem \ref{t2}.
\end{proof}

\subsection{Proof of Theorem \ref{t3}}

\begin{proof}
The proof of (ii) is obvious. On the other hand, the proof of (i) is very
similar to the proof of clause (i) of Theorem \ref{t2}, so we skip its
beginning, and state the starting system of equations, obtained with the same
choice of variables as in equations (\ref{e1})-(\ref{e4}):
\begin{align*}
\lambda x &  =ky,\\
\lambda y &  =\left(  k+1\right)  x+\left(  k-1\right)  y+\left(
n-2k-3\right)  z+2t,\\
\lambda z &  =ky+\left(  n-2k-4\right)  z+2t,\\
\lambda t &  =ky+\left(  n-2k-3\right)  z.
\end{align*}
Solving this system with respect to $y,$ we find that%
\begin{align}
x &  =\frac{k}{\lambda}y,\label{eq1}\\
z &  =\left(  1-\frac{k\left(  k+1\right)  }{\lambda\left(  \lambda+1\right)
}\right)  y,\label{eq2}\\
t &  =\frac{\lambda+1}{\lambda+2}\left(  1-\frac{k\left(  k+1\right)
}{\lambda\left(  \lambda+1\right)  }\right)  y.\label{eq3}%
\end{align}
Assume for a contradiction that $\lambda\geq n-k-2$. Proceeding further as in
the proof of clause (i) of Theorem \ref{t2}, we get the inequalities
\begin{align*}
\left\langle A\left(  K_{n-k-1}\right)  \mathbf{x}^{\prime\prime}%
,\mathbf{x}^{\prime\prime}\right\rangle  &  =\left\langle A\left(  G\right)
\mathbf{x},\mathbf{x}\right\rangle +2t^{2}-2k\left(  k+1\right)
xy=\lambda+2t^{2}-2k\left(  k+1\right)  xy\\
\left\langle A\left(  K_{n-k-1}\right)  \mathbf{x}^{\prime\prime}%
,\mathbf{x}^{\prime\prime}\right\rangle  &  \leq\left(  n-k-2\right)  \left(
1-\left(  k+1\right)  x^{2}\right)  \leq\lambda\left(  1-\left(  k+1\right)
x^{2}\right)  .
\end{align*}
Hence,%
\[
2k\left(  k+1\right)  xy-\left(  k+1\right)  \lambda x^{2}\geq2t^{2}.
\]
Using (\ref{eq1}), (\ref{eq2}) and (\ref{eq3}), after simple algebra we get
\[
k^{2}\left(  k+1\right)  >2\lambda\left(  1-\frac{1}{\lambda+2}\right)
^{2}\left(  1-\frac{k\left(  k+1\right)  }{\lambda\left(  \lambda+1\right)
}\right)  ^{2}.
\]
Applying the Bernoulli inequality to the right-side of the above inequality,
we find that
\[
k^{2}\left(  k+1\right)  >2\lambda-\frac{4\lambda}{\lambda+2}-\frac{4k\left(
k+1\right)  }{\lambda+1}\geq2\lambda-\frac{4\lambda+4k\left(  k+1\right)
}{\lambda+1}%
\]
Now in view of \ $\lambda\geq n-k-2\geq k^{3}/2+k^{2}/2+3,$ we easily see
that
\[
\frac{4\lambda+4k\left(  k+1\right)  }{\lambda+1}\leq6,
\]
and so,%
\[
k^{2}\left(  k+1\right)  >2\lambda-6\geq k^{3}+k^{2},
\]
a contradiction, completing the proof of Theorem \ref{t3}.\medskip
\end{proof}

\subsection{Proof of Theorem \ref{mtc}}

\begin{proof}
Let $k\geq1,$ $n\geq k^{3}/4+k+4,$ and let $G$ be a graph of order $n$, with
$\delta\left(  G\right)  \geq k.$ Write $m$ for the number of edges of $G,$
and set $\delta:=\delta\left(  G\right)  .$

Assume that $\lambda\left(  G\right)  \geq n-k-1,$ but $G$ has no Hamiltonian
cycle. To prove the theorem we need to show that $G=L_{k}\left(  n\right)  $
or $G=M_{k}\left(  n\right)  .$ Note that, in view of Theorem \ref{t2}, it is
sufficient to prove that $\mathrm{cl}_{n}\left(  G\right)  =$ $L_{k}\left(
n\right)  $ or $\mathrm{cl}_{n}\left(  G\right)  =M_{k}\left(  n\right)  ,$ so
this will be our main goal.

Clearly, $\mathrm{cl}_{n}\left(  G\right)  $\ has no Hamiltonian cycle and%
\[
\delta\left(  \mathrm{cl}_{n}\left(  G\right)  \right)  \geq\delta\left(
G\right)  \geq k,\text{ \ \ \ \ }\lambda\left(  \mathrm{cl}_{n}\left(
G\right)  \right)  \geq\lambda\left(  G\right)  \geq n-k-1,
\]
so for the rest of the proof we assume that $G=\mathrm{cl}_{n}\left(
G\right)  .$ The main consequence of this assumption is that
\begin{equation}
d_{i}+d_{j}\leq n-1 \label{clin}%
\end{equation}
for every two nonadjacent vertices $i$ and $j$.

Next, since $G$ has no Hamiltonian cycle, Theorem \ref{cht} implies that there
is an integer $s<n/2$ such that $d_{s}\leq s$ and $d_{n-s}\leq n-s-1.$
Obviously, $s\geq\delta\geq k,$ and we easily find an upper bound on $2m:$
\begin{align*}
2m  &  =\sum_{i=1}^{s}d_{i}+\sum_{i=s+1}^{n-s}d_{i}+\sum_{i=n-s+1}^{n}d_{i}\\
&  \leq s^{2}+\left(  n-2s\right)  \left(  n-s-1\right)  +s\left(  n-1\right)
\\
&  =n^{2}-2sn+3s^{2}+s-n.
\end{align*}
Clearly, the expression $n^{2}-2sn+3s^{2}+s-n$ is convex in $s$; hence it is
maximal in $s$ for $s=\delta$ or $s=\left(  n-1\right)  /2.$ Hence, either
\begin{equation}
2m\leq n^{2}-2\delta n+3\delta^{2}+\delta-n \label{u1}%
\end{equation}
or
\begin{equation}
2m\leq n^{2}-\left(  n-1\right)  n+3\frac{\left(  n-1\right)  ^{2}}{4}%
+\frac{\left(  n-1\right)  }{2}-n. \label{u2}%
\end{equation}

On the other hand, inequality\ (\ref{in}) implies that
\[
n-k-1\leq\lambda\left(  G\right)  \leq\frac{\delta-1}{2}+\sqrt{2m-n\delta
+\frac{\left(  \delta+1\right)  ^{2}}{4}}.
\]
Hence, in view of Proposition \ref{pro}, we get
\[
n-k-1\leq\lambda\left(  G\right)  \leq\frac{k-1}{2}+\sqrt{2m-nk+\frac{\left(
k+1\right)  ^{2}}{4}},\text{ }%
\]
which, after some algebra, gives%
\begin{equation}
2m\geq n^{2}-2kn+2k^{2}+k-n.\label{lb}%
\end{equation}

Next, we prove that $s=k.$ Indeed, if $s\geq k+1,$ then (\ref{u1}) and
(\ref{u2}) imply that either
\[
n^{2}-2\left(  k+1\right)  n+3\left(  k+1\right)  ^{2}+k+1-n\geq
n^{2}-2kn+2k^{2}+k-n
\]
or%
\[
n^{2}-\left(  n-1\right)  n+\frac{3}{4}\left(  n-1\right)  ^{2}+\frac{1}%
{2}\left(  n-1\right)  -n\geq n^{2}-2kn+2k^{2}+k-n.
\]
Each of these inequalities leads to a contradiction, so we have $s=k,$ and
thus $\delta=k.$ Therefore,%
\[
d_{1}=\cdots=d_{k}=k.
\]
\qquad

Our next goal is to show that $d_{k+1}\geq n-k-1-k^{2}.$ \ Indeed, suppose
that
\[
d_{k+1}<n-k-1-k^{2}.
\]
Now, using Theorem \ref{cht}, we get%
\begin{align*}
2m &  =\sum_{i=1}^{k}d_{i}+d_{k+1}+\sum_{i=k+2}^{n-k}d_{i}+\sum_{i=n-k+1}%
^{n}d_{i}\\
&  <k^{2}+n-k-1-k^{2}+\left(  n-2k-1\right)  \left(  n-k-1\right)  +k\left(
n-1\right)  \\
&  =n^{2}-2kn+2k^{2}+k-n,
\end{align*}
contradicting (\ref{lb}). Hence $d_{i}\geq n-k-1-k^{2}$ for every
$i\in\left\{  k+1,\ldots,n\right\}  $.

Next, we shall show that the vertices $k+1,\ldots,n$ induce a complete graph
in $G.$ Indeed, let $i\in\left\{  k+1,\ldots,n\right\}  $ and $j\in\left\{
k+1,\ldots,n\right\}  $ be two distinct vertices of $G.$ If they are
nonadjacent, then
\begin{align*}
d_{i}+d_{j} &  \geq2n-2k-2-2k^{2}\\
&  \geq n+k^{3}+k+4-2k-2-2k^{2}\\
&  =n-1+k^{3}-2k^{2}-k+3>n-1,
\end{align*}
contradicting (\ref{clin}).

Write $X$ for the vertex set $\left\{  1,\ldots,k\right\}  .$ Write $Y$ for
the set of vertices in $\left\{  k+1,\ldots n\right\}  $ having neighbors in
$X.$ It is easy to see that $Y\neq\varnothing,$ since $\left\vert X\right\vert
=k$ and so any vertex in $X$ must have a neighbor in $\left\{  k+1,\ldots
,n\right\}  .$

In fact, every vertex from $Y$ is adjacent to every vertex in $X$. Indeed,
suppose that this is not the case, and let $w\in\left\{  k+1,\ldots,n\right\}
,$ $u\in X,$ $v\in X$ be such that $w$ is adjacent to $u,$ but not to $v.$ We
see that
\[
d_{w}+d_{v}\geq n-k+k=n,
\]
contradicting (\ref{clin}).

Next, let $l:=\left\vert Y\right\vert $ and note that $1\leq l\leq k,$ since
$d_{1}=k.$ If $l=1,$ then $G=L_{k}\left(  n\right)  ,$ and if $l=k,$ then
$G=M_{k}\left(  n\right)  .$ To finish the proof we shall show that if
$1<l<k,$ then $G$ has a Hamiltonian cycle, which contradicts the assumptions
about $G$.

Indeed, let $H$ be the graph induced by the set $X\cup Y.$ Since $K_{l}%
\vee\overline{K}_{k}\subset H$ and $l\geq2,$ we see that $H$ is $2$-connected.
Further, if $u$ and $v$ are distinct nonadjacent vertices of $H,$ with degrees
$d_{u}^{\prime}$ and $d_{v}^{\prime},$ they must belong to $X,$ and so
$d_{u}^{\prime}=d_{u}=k$ and $d_{v}^{\prime}=d_{v}=k.$ That is to say,
\[
d_{u}^{\prime}+d_{v}^{\prime}=2k>k+l.
\]
Theorem \ref{ot} implies that $H$ is Hamiltonian-connected, and it is easy to
see that $G$ has a Hamiltonian cycle. The proof of Theorem \ref{mtc} is completed.
\end{proof}

\subsection{Proof of Theorem \ref{mtp}}

\begin{proof}
Although this proof is very close to the proof of Theorem \ref{mtc}, we shall
carry it in full, due to the numerous specific details.

Let $k\geq1,$ $n\geq k^{3}/2+k^{2}/2+k+5$ and let $G$ be a graph of order $n$,
with $\delta\left(  G\right)  \geq k.$ Write $m$ for the number of edges of
$G$, and set $\delta:=\delta\left(  G\right)  .$

Assume that $\lambda\left(  G\right)  \geq n-k-2,$ but $G$ has no Hamiltonian
path. To prove the theorem we need to show that $G=N_{k}\left(  n\right)  $ or
$G=K_{n-k-1}+K_{k+1}.$ Note that, in view of Theorem \ref{t2}, it is
sufficient to prove that $\mathrm{cl}_{n}\left(  G\right)  =$ $N_{k}\left(
n\right)  $ or $\mathrm{cl}_{n}\left(  G\right)  =K_{n-k-1}+K_{k+1},$ so this
will be our main goal.

Clearly, $\mathrm{cl}_{n}\left(  G\right)  $\ has no Hamiltonian path and%
\[
\delta\left(  \mathrm{cl}_{n}\left(  G\right)  \right)  \geq\delta\left(
G\right)  \geq k,\text{ \ \ \ }\lambda\left(  \mathrm{cl}_{n}\left(  G\right)
\right)  \geq\lambda\left(  G\right)  \geq n-k-2,
\]
so for the rest of the proof we assume that $G=\mathrm{cl}_{n}\left(
G\right)  .$ The main consequence of this assumption is that
\begin{equation}
d_{i}+d_{j}\leq n-2 \label{clinp}%
\end{equation}
for every two nonadjacent vertices $i$ and $j$.

Next, since $G$ has no Hamiltonian path, Corollary \ref{chc} implies that
there is an integer $s\leq n/2$ such that $d_{s}\leq s-1$ and $d_{n-s+1}\leq
n-s.$ Obviously, $s\geq\delta+1\geq k+1,$ and we easily find an upper bound on
$2m:$
\begin{align*}
2m  &  =\sum_{i=1}^{s}d_{i}+\sum_{i=s+1}^{n-s+1}d_{i}+\sum_{i=n-s+2}^{n}%
d_{i}\\
&  \leq s\left(  s-1\right)  +\left(  n-2s+1\right)  \left(  n-s-1\right)
+\left(  s-1\right)  \left(  n-1\right) \\
&  =n^{2}-2sn+3s^{2}-s-n.
\end{align*}
Clearly, the expression $n^{2}-2sn+3s^{2}-s-n$ is convex in $s$; hence it is
maximal in $s$ for $s=\delta+1$ or $s=n/2.$ Hence, either
\begin{equation}
2m\leq n^{2}-2\left(  \delta+1\right)  n+3\left(  \delta+1\right)
^{2}-\left(  \delta+1\right)  -n \label{u21}%
\end{equation}
or
\begin{equation}
2m\leq\frac{3}{4}n^{2}-\frac{3}{2}n. \label{u22}%
\end{equation}

On the other hand, inequality\ (\ref{in}) implies that
\[
n-k-2\leq\lambda\left(  G\right)  \leq\frac{\delta-1}{2}+\sqrt{2m-n\delta
+\frac{\left(  \delta+1\right)  ^{2}}{4}}.
\]
In view of Proposition \ref{pro}, we get
\[
n-k-2\leq\lambda\left(  G\right)  \leq\frac{k-1}{2}+\sqrt{2m-nk+\frac{\left(
k+1\right)  ^{2}}{4}},\text{ }%
\]
which, after some algebra, gives%
\begin{equation}
2m\geq n^{2}-2kn+2k^{2}+4k-3n+2. \label{lb1}%
\end{equation}

We shall prove that $s=k+1.$ Indeed, if $s\geq k+2,$ then (\ref{u21}) and
(\ref{u22}) imply that either
\[
n^{2}-2\left(  k+2\right)  n+3\left(  k+2\right)  ^{2}-\left(  k+2\right)
-n\geq n^{2}-2kn+2k^{2}+4k-3n+2
\]
or%
\[
\frac{3}{4}n^{2}-\frac{3}{2}n\geq n^{2}-2kn+2k^{2}+4k-3n+2.
\]
Both of these inequalities lead to a contradiction, so we have $s=k+1$ and
thus $\delta=k.$ Therefore,%
\[
d_{1}=\cdots=d_{k+1}=k.
\]

Our next goal is to show that $d_{k+2}\geq n-2k-2-k^{2}.$ \ Indeed, suppose
that
\[
d_{k+2}<n-2k-2-k^{2}.
\]
Now we get%
\begin{align*}
2m  &  =\sum_{i=1}^{k+1}d_{i}+d_{k+2}+\sum_{i=k+3}^{n-k}d_{i}+\sum
_{i=n-k+1}^{n}d_{i}\\
&  <k\left(  k+1\right)  +n-2k-2-k^{2}+\left(  n-2k-2\right)  \left(
n-k-2\right)  +k\left(  n-1\right) \\
&  =n^{2}-2kn+2k^{2}+4k-3n+2,
\end{align*}
contradicting (\ref{lb1}). Hence $d_{i}\geq n-2k-2-k^{2}$ for every
$i\in\left\{  k+2,\ldots,n\right\}  $.

Next, we show that the vertices $k+2,\ldots,n$ induce a complete graph in $G.$
Indeed, let $i\in\left\{  k+2,\ldots,n\right\}  $ and $j\in\left\{
k+2,\ldots,n\right\}  $ be two distinct vertices of $G.$ If they are
nonadjacent, then
\begin{align*}
d_{i}+d_{j} &  \geq2n-4k-4-2k^{2}\\
&  \geq n+k^{3}+k^{2}+2k+5-4k-4-2k^{2}\\
&  =\left(  n-2\right)  +k^{3}-k^{2}-2k+3>n-2,
\end{align*}
contradicting (\ref{clinp}).

Write $X$ for the vertex set $\left\{  1,\ldots,k+1\right\}  .$ Write $Y$
for\ the set of vertices in $\left\{  k+2,\ldots n\right\}  $ that have
neighbors in $X.$ If $Y=\varnothing,$ then $G$ is a disconnected graph and the
order of its largest component is at most $n-k-1$. Also $X$ induces a
$K_{k+1}.$ Clearly, the inequality $\lambda\left(  G\right)  \geq n-k-2$
implies that $G=K_{n-k-1}+K_{k+1},$ completing the proof if $Y=\varnothing$.

Now, suppose that $Y\neq\varnothing$ $.$ We shall show that every vertex in
$Y$ is adjacent to every vertex in $X$. Indeed, suppose that this is not the
case, and let $w\in\left\{  k+1,\ldots,n\right\}  ,$ $u\in X,$ $v\in X$ be
such that $w$ is adjacent to $u,$ but not to $v.$ We see that
\[
d_{w}+d_{v}\geq n-k-1+k=n-1,
\]
contradicting (\ref{clin}).

Next, let $l:=\left\vert Y\right\vert $ and note that $1\leq l\leq k,$ as
$d_{1}=k.$ If $l=k,$ then $G=N_{k}\left(  n\right)  .$ To finish the proof we
shall show that if $1\leq l<k,$ then $G$ has a Hamiltonian path.

Indeed, let $H$ be the graph induced by the set $X\cup Y.$ Further, if $u$ and
$v$ are distinct nonadjacent vertices of $H,$ with degrees $d_{u}^{\prime}$
and $d_{v}^{\prime},$ they must belong to $X,$ and so $d_{u}^{\prime}=d_{u}=k$
and $d_{v}^{\prime}=d_{v}=k.$ That is to say,
\[
d_{u}^{\prime}+d_{v}^{\prime}=2k\geq k+1+l.
\]
Theorem \ref{ot1} implies that $H$ contains a Hamiltonian cycle, and hence $G$
has a Hamiltonian path. The proof of Theorem \ref{mtp} is completed.
\end{proof}

\subsection{Proofs of Propositions \ref{pro1} and \ref{pro2}}

\begin{proof}
[\textbf{Proof of Proposition \ref{pro1}}]Suppose that $2k+1\leq n\leq
k^{3}/2+k+1$, and consider the graph $M_{k}\left(  n\right)  $. Write $X$ for
the set of vertices of $M_{k}\left(  n\right)  $ of degree $k,$ let $Y$ be the
set of their neighbors, and let $Z$ be the set of the remaining $n-2k$
vertices of $M_{k}\left(  n\right)  .$ Select two distinct vertices $u\in
Z\cup Y$ and $v\in Z\cup Y,$ remove the edge $\left\{  u,v\right\}  $, and
write $G$ for the resulting graph.

To begin with, note that $\delta\left(  G\right)  \geq k$ and $G$ is a proper
subgraph of $M_{k}\left(  n\right)  .$ To complete the proof we define a unit
$n$-vector $\mathbf{x}$ and show that
\[
\lambda\left(  G\right)  >\left\langle A\left(  G\right)  \mathbf{x}%
,\mathbf{x}\right\rangle \geq n-k-1.
\]
Thus, let
\[
x=\frac{2}{k^{2}\sqrt{n-k+4/k^{3}}},\text{ \ \ }y=\frac{1}{\sqrt{n-k+4/k^{3}}%
},
\]
and define the $n$-vector $\mathbf{x}:=$ $\left(  x_{1},\ldots,x_{n}\right)  $
as%
\[
x_{i}:=\left\{
\begin{array}
[c]{ll}%
x\text{,} & \text{if }i\in X\text{;}\\
y\text{,} & \text{if }i\in Y\cup Z\text{.}%
\end{array}
\right.
\]
Note that $\left\Vert \mathbf{x}\right\Vert =1,$ because
\[
\left\Vert \mathbf{x}\right\Vert ^{2}=kx^{2}+\left(  n-k\right)  y^{2}%
=\frac{4}{k^{3}\left(  n-k+4/k^{3}\right)  }+\frac{n-k}{n-k+4/k^{3}}=1.
\]
Therefore,%
\begin{align*}
\lambda\left(  G\right)   &  \geq\left\langle A\left(  G\right)
\mathbf{x},\mathbf{x}\right\rangle =\left(  n-k\right)  \left(  n-k-1\right)
y^{2}-2y^{2}+2k^{2}xy\\
&  =\left(  n-k\right)  \left(  n-k-1\right)  y^{2}+2y^{2}\\
&  =\frac{n-k}{n-k+4/k^{3}}\left(  n-k-1\right)  +\frac{2}{n-k+4/k^{3}}\\
&  =n-k-1-\frac{4}{k^{3}\left(  n-k+4/k^{3}\right)  }\left(  n-k-1-\frac
{k^{3}}{2}\right)  \\
&  \geq n-k-1.
\end{align*}
Finally, the inequality $\lambda\left(  G\right)  \geq\left\langle A\left(
G\right)  \mathbf{x},\mathbf{x}\right\rangle $ is strict, because obviously
$\mathbf{x}$ is not an eigenvector to $\lambda\left(  G\right)  .$
\end{proof}

$\medskip$

\begin{proof}
[\textbf{Proof of Proposition \ref{pro2}}]Suppose that $2k+1\leq n\leq
k^{3}/2+k^{2}/2+k+2$, and consider the graph $N_{k}\left(  n\right)  .$ Write
$X$ for the set of vertices of $N_{k}\left(  n\right)  $ of degree $k,$ let
$Y$ be the set of their neighbors, and let $Z$ be the set of the remaining
$n-2k-1$ vertices of $N_{k}\left(  n\right)  .$ Select two distinct vertices
$u\in Z\cup Y$ and $v\in Z\cup Y,$ remove the edge $\left\{  u,v\right\}  $,
and write $G$ for the resulting graph.

To begin with, note that $\delta\left(  G\right)  \geq k$ and $G$ is a proper
subgraph of $N_{k}\left(  n\right)  .$ To complete the proof we define a unit
$n$-vector $\mathbf{x}$\ and show that%
\[
\lambda\left(  G\right)  >\left\langle A\left(  G\right)  \mathbf{x}%
,\mathbf{x}\right\rangle \geq n-k-2.
\]
Thus, let%
\[
x=\frac{2}{k\left(  k+1\right)  \sqrt{n-k-1+4/k^{2}\left(  k+1\right)  }%
},\ \ \text{\ }y=\frac{1}{\sqrt{n-k-1+4/k^{2}\left(  k+1\right)  }},
\]
and define the $n$-vector $\mathbf{x}:=$ $\left(  x_{1},\ldots,x_{n}\right)  $
as%
\[
x_{i}:=\left\{
\begin{array}
[c]{ll}%
x\text{,} & \text{if }i\in X\text{;}\\
y\text{,} & \text{if }i\in Y\cup Z\text{.}%
\end{array}
\right.
\]
Note that $\left\Vert \mathbf{x}\right\Vert =1,$ because
\begin{align*}
\left\Vert \mathbf{x}\right\Vert ^{2} &  =\left(  k+1\right)  x^{2}+\left(
n-k-1\right)  y^{2}\\
&  =\frac{4}{k^{2}\left(  k+1\right)  \left(  n-k-1+4/k^{2}\left(  k+1\right)
\right)  }+\frac{n-k-1}{n-k-1+4/k^{2}\left(  k+1\right)  }=1.
\end{align*}
Therefore,%
\begin{align*}
\lambda\left(  G\right)   &  \geq\left\langle A\left(  G\right)
\mathbf{x},\mathbf{x}\right\rangle =\left(  n-k-1\right)  \left(
n-k-2\right)  y^{2}-2y^{2}+2k\left(  k+1\right)  xy\\
&  =\left(  n-k-1\right)  \left(  n-k-2\right)  y^{2}+2y^{2}\\
&  =\frac{\left(  n-k-1\right)  \left(  n-k-2\right)  }{n-k-1+4/k^{2}\left(
k+1\right)  }+\frac{2}{n-k-1+4/k^{2}\left(  k+1\right)  }\\
&  =n-k-2-\frac{4}{k^{2}\left(  k+1\right)  \left(  n-k-1+4/k^{2}\left(
k+1\right)  \right)  }\left(  n-k-2-\frac{k^{2}\left(  k+1\right)  }%
{2}\right)  \\
&  \geq n-k-2.
\end{align*}
Finally, the inequality $\lambda\left(  G\right)  \geq\left\langle A\left(
G\right)  \mathbf{x},\mathbf{x}\right\rangle $ is strict, because obviously
$\mathbf{x}$ is not an eigenvector to $\lambda\left(  G\right)  .$
\end{proof}

\section{\label{s4}Concluding remarks}

It should be noted that most of the results discussed in the present paper and
in the references \cite{Ben15,GeNi16, Li15, LiNi15, LLT12, LSX15, NiGe15,
Zho10} deal exclusively with very dense graphs, which makes these results
somewhat one-sided. Hoping to change this tendency, we would like to state two
open problems.\medskip

First, recall that Dirac's theorem \cite{Dir52} is probably the most famous
sufficient condition for Hamiltonian cycles. Yet, no comparable spectral
statement seems to be known so far.

\begin{problem}
Find a spectral sufficient condition for Hamiltonian cycles that would imply
Dirac's sufficient condition.
\end{problem}

Second, a deep result of Krivelevich and Sudakov \cite{KrSu03} establishes a
sufficient condition on the second largest singular value of a regular graph
that implies existence of Hamiltonian cycles. Two attempts have been made to
extend this result to nonregular graphs, but these extensions forsake the
adjacency matrix for other matrices (\cite{BuCh10}, \cite{FaYu12}), so
comparisons are difficult.

Hence, it is worth to reiterate the following problem, first raised in
\cite{KrSu03}:

\begin{problem}
Extend the result of Krivelevich and Sudakov to nonregular graphs, using the
second largest singular value of the adjacency matrix.
\end{problem}

\end{document}